%%%%%%%%%%%%%%%%%%%%%%%%%%%%%%%%%%%%%%%%%%%%%%%%%%%%%%%%%%
%%%%%%%%%%%%%%%%%%%%%%%%%%%%%%%%%%%%%%%%%%%%%%%%%%%%%%%%%%
%%
%%     This is the AMS-LaTeX file:
%%
%%     Colli-Gilardi-Sprekels 22
%%     Longtime behavior for a generalized Cahn--Hilliard system with fractional operators
%%     
%%
%%%%%%%%%%%%%%%%%%%%%%%%%%%%%%%%%%%%%%%%%%%%%%%%%%%%%%%%%%

\def\LaTeX{%
  \let\Begin\begin
  \let\End\end
  \let\salta\relax
  \let\finqui\relax
  \let\futuro\relax}

\def\UK{\def\our{our}\let\sz s}
\def\USA{\def\our{or}\let\sz z}

\UK
%\USA

%%%%%%%%%%%%%%%%%%%%%%%%%%%%%%%%%

% scegliere fra \TeX e \LaTeX  e fra  \UK oppure \USA

%\TeX
\LaTeX

%\UK
\USA

%%%%%%%%%%%%%%%%%%%%%%%%%%%%%%%%%
%% page layout
%%%%%%%%%%%%%%%%%%%%%%%%%%%%%%%%%

\salta

\documentclass[twoside,12pt]{article}
\setlength{\textheight}{24cm}
\setlength{\textwidth}{16cm}
\setlength{\oddsidemargin}{2mm}
\setlength{\evensidemargin}{2mm}
\setlength{\topmargin}{-15mm}
\parskip2mm

%%%%%%%%%%%%%%%%%%%%%%%%%%%%%%%%%
%% packages
%%%%%%%%%%%%%%%%%%%%%%%%%%%%%%%%%

%\usepackage{color}
\usepackage[usenames,dvipsnames]{color}
\usepackage{amsmath}
\usepackage{amsthm}
\usepackage{amssymb}
\usepackage[mathcal]{euscript}
\usepackage{cite}

%\usepackage[notref,notcite]{showkeys}
%\usepackage{showkeys}
%
%		COLORS FOR CORRECTIONS
%
% do the same, please (i.e., don't use the standard {\color{red} text} or similar): 
% just choose the color you prefer in \def\yourname

% example of use:  \juerg{I want this to become blue}

\definecolor{viola}{rgb}{0.3,0,0.7}
\definecolor{ciclamino}{rgb}{0.5,0,0.5}

\def\pier #1{{\color{red}#1}}
\def\juerg #1{{\color{blue}#1}}

\def\pier #1{#1}
\def\juerg #1{#1}

%\def\revis #1{{\color{red}#1}}
%\def\revis #1{#1}

%%%%%%%%%%%%%%%%%%%%%%%%%%%%%%%%%
%% you may adjust the baseline
%%%%%%%%%%%%%%%%%%%%%%%%%%%%%%%%%

%%%%%%%%%%%%%%%%%%%%%%%%%%%%%%%%%
%% bibliographystyle
%%%%%%%%%%%%%%%%%%%%%%%%%%%%%%%%%

\bibliographystyle{plain}

%%%%%%%%%%%%%%%%%%%%%%%%%%%%%%%%%
%% environments
%%%%%%%%%%%%%%%%%%%%%%%%%%%%%%%%%

%

\finqui

\def\Beq{\Begin{equation}}
\def\Eeq{\End{equation}}
\def\Bsist{\Begin{eqnarray}}
\def\Esist{\End{eqnarray}}

\def\Bthm{\Begin{theorem}}
\def\Ethm{\End{theorem}}

\def\Bprop{\Begin{proposition}}
\def\Eprop{\End{proposition}}

\def\Bex{\Begin{example}\rm}
\def\Eex{\End{example}}
\def\Bnot{\Begin{notation}\rm}
\def\Enot{\End{notation}}

\def\Bcenter{\Begin{center}}
\def\Ecenter{\End{center}}
\let\non\nonumber

%%%%%%%%%%%%%%%%%%%%%%%%%%%%%%%%%
%% macros
%%%%%%%%%%%%%%%%%%%%%%%%%%%%%%%%%

% macro salvate

% sottosezioni non numerate

\def\step #1 \par{\medskip\noindent{\bf #1.}\quad}

% abbreviazioni di parole

\def\Lip{Lip\-schitz}

\def\aand{\quad\hbox{and}\quad}

\def\rhs{right-hand side}

\def\omegalimit{$\omega$-limit}

% versioni inglesi (UK) o americane (USA)

\def\bhv{behavi\our}

% bold, cal e mathop

\def\multibold #1{\def\arg{#1}%
  \ifx\arg\pto \let\next\relax
  \else
  \def\next{\expandafter
    \def\csname #1#1#1\endcsname{{\bf #1}}%
    \multibold}%
  \fi \next}

\def\pto{.}

\def\multical #1{\def\arg{#1}%
  \ifx\arg\pto \let\next\relax
  \else
  \def\next{\expandafter
    \def\csname cal#1\endcsname{{\cal #1}}%
    \multical}%
  \fi \next}

% operatori

\def\multimathop #1 {\def\arg{#1}%
  \ifx\arg\pto \let\next\relax
  \else
  \def\next{\expandafter
    \def\csname #1\endcsname{\mathop{\rm #1}\nolimits}%
    \multimathop}%
  \fi \next}

\multibold
qwertyuiopasdfghjklzxcvbnmQWERTYUIOPASDFGHJKLZXCVBNM.

\multical
QWERTYUIOPASDFGHJKLZXCVBNM.

\multimathop
diag dist div dom mean meas sign supp .

% accorpamenti di formule citate:
% uso  \accorpa {prima}{seconda}
%      \Accorpa\cs prima seconda (con il comodo blank anche dopo)
% NB: \Accorpa definisce \cs come l'accorpamento delle due citazioni
% e scrive sul file.log

\def\accorpa #1#2{\eqref{#1}--\eqref{#2}}
\def\Accorpa #1#2 #3 {\gdef #1{\eqref{#2}--\eqref{#3}}%
  \wlog{}\wlog{\string #1 -> #2 - #3}\wlog{}}

% macro comode

\def\separa{\noalign{\allowbreak}}

\def\somma #1#2#3{\sum_{#1=#2}^{#3}}

\def\graffe #1{\mathopen\{#1\mathclose\}}

\def\<#1>{\mathopen\langle #1\mathclose\rangle}
\def\norma #1{\mathopen \| #1\mathclose \|}
\def\Norma #1{\Bigl\| #1 \Bigr\|}

\def\[#1]{\mathopen\langle\!\langle #1\mathclose\rangle\!\rangle}

\def\ioT {\int_0^T}

\def\intQ{\int_{\Omega\times(0,T)}\hskip-.5em}
\def\iO{\int_\Omega}

\def\dt{\partial_t}

\def\cpto{\,\cdot\,}

\def\checkmmode #1{\relax\ifmmode\hbox{#1}\else{#1}\fi}
\def\aeO{\checkmmode{a.e.\ in~$\Omega$}}

\def\aat{\checkmmode{for a.a.~$t\in(0,T)$}}

\def\Aat{\checkmmode{for a.a.~$t\in(0,+\infty)$}}

% insiemi numerici

\def\erre{{\mathbb{R}}}

% spazi di funzioni a valori vettoriali su [0,T], [0,t], [0,s], [0,+\infty), [\delta,T]

% Come ricordare: in generale i simboli L H W  C da soli per gli spazi su (0,T)
% gli stessi raddoppiati per (0,+\infty)
% aggiunta di t o s al simbolo per (0,t) e (0,s)
% aggiunta di d al simbolo semplice o doppio per intervalli (\delta,T) e (\delta,+\infty)
% il simbolo C e i suoi derivati mettono le quadre anziche' le tonde

% Esempi   \L2V   \L\infty\Vp   \W{1,1}H   \C0H   \LL2V   \CC0\Vp   \Ld2V  \CCdH

\def\genspazio #1#2#3#4#5{#1^{#2}(#5,#4;#3)}
\def\spazio #1#2#3{\genspazio {#1}{#2}{#3}T0}

\def\L {\spazio L}
\def\H {\spazio H}
\def\W {\spazio W}

\def\C #1#2{C^{#1}([0,T];#2)}
\def\spazioinf #1#2#3{\genspazio {#1}{#2}{#3}{+\infty}0}
\def\LL {\spazioinf L}

% spazi di funzioni su \Omega, \Gamma, Q e \Sigma

\def\Lx #1{L^{#1}(\Omega)}

\def\Luno{\Lx 1}
\def\Ldue{\Lx 2}
\def\Linfty{\Lx\infty}

% spazi di funzioni su Q e S

% lettere greche

\let\theta\vartheta

\let\phi\varphi

\let\hat\widehat

\let\TeXchi\chi                         % new \chi, exactly on the baseline
\newbox\chibox
\setbox0 \hbox{\mathsurround0pt $\TeXchi$}
\setbox\chibox \hbox{\raise\dp0 \box 0 }
\def\chi{\copy\chibox}

% quadratino di fine dimostrazione

\def\QED{\hfill $\square$}

% abbreviazioni specifiche del lavoro

\def\VA #1{V_A^{#1}}
\def\VB #1{V_B^{#1}}

\def\Beta{\hat\beta}
\def\betal{\beta_\lambda}
\def\Betal{\hat\beta_\lambda}

\def\Pi{\hat\pi}
\def\Lpi{L_\pi}

\def\yl{y^\lambda}
\def\mul{\mu^\lambda}

\def\yh{\hat y_h}

\def\overyh{\overline y_h}
\def\overmuh{\overline\mu_h}
\def\underyh{\underline y_h}
\def\undermuh{\underline\mu_h}

\def\yn{y^n}
\def\ynp{y^{n+1}}

\def\mun{\mu^n}
\def\munp{\mu^{n+1}}

\def\un{u^n}
\def\unp{u^{n+1}}
\def\dhyn{\frac{\ynp-\yn}h}

\def\yk{y^k}
\def\muk{\mu^k}
\def\uk{u^k}

\def\yz{y_0}
\def\mz{m_0}

\def\Lpi{L_\pi}
\def\yu{y^1}
\def\muu{\mu^1}

\def\yo{y_\omega}
\def\tn{t_n}
\def\yi{y_\infty}
\def\mui{\mu_\infty}
\def\ui{u_\infty}

%%%%%%%%%%%%%%%%%%%%%%%%%%%%%%
\Begin{document}
%%%%%%%%%%%%%%%%%%%%%%%%%%%%%%%%%

%%%%%%%%%%%%%%%%%%%%%%%%%%%%%%%%%
%% front page
%%%%%%%%%%%%%%%%%%%%%%%%%%%%%%%%%

%
\title{Longtime behavior
  \\ for a generalized Cahn--Hilliard system
  \\ with fractional operators}
\author{}
\date{}
\maketitle
\Bcenter
\vskip-1.5cm
{\large\sc Pierluigi Colli$^{(1)}$}\\
{\normalsize e-mail: {\tt pierluigi.colli@unipv.it}}\\[.25cm]
{\large\sc Gianni Gilardi$^{(1)}$}\\
{\normalsize e-mail: {\tt gianni.gilardi@unipv.it}}\\[.25cm]
{\large\sc J\"urgen Sprekels$^{(2)}$}\\
{\normalsize e-mail: {\tt sprekels@wias-berlin.de}}\\[.45cm]
$^{(1)}$
{\small Dipartimento di Matematica ``F. Casorati'', Universit\`a di Pavia}\\
{\small and Research Associate at the IMATI -- C.N.R. Pavia}\\
{\small via Ferrata 5, 27100 Pavia, Italy}\\[.2cm]
$^{(2)}$
{\small Department of Mathematics}\\
{\small Humboldt-Universit\"at zu Berlin}\\
{\small Unter den Linden 6, 10099 Berlin, Germany}\\[2mm]
{\small and}\\[2mm]
{\small Weierstrass Institute for Applied Analysis and Stochastics}\\
{\small Mohrenstrasse 39, 10117 Berlin, Germany}
\Ecenter
\begin{center}
\emph{Dedicated to Professor Antonino Maugeri \\on the occasion of his 75th birthday\\ with best wishes}
\end{center}
\Begin{abstract}\noindent
In this contribution, 
we \pier{deal with} the longtime \bhv\ of the solutions
to the fractional variant of the Cahn--Hilliard system\pier{, 
with possibly singular potentials, that we have 
recently investigated in the paper
{\sl Well-posedness and regularity for a generalized 
fractional Cahn--Hilliard system\/}.}
More precisely, we study the \omegalimit\ of the phase parameter $y$
and characterize it completely.
Our characterization depends on the first eigenvalues $\lambda_1\geq0$ 
of one of the operators involved:
if $\lambda_1>0$, then the chemical potential $\mu$ vanishes at infinity
and every element $\yo$ of the \omegalimit\ is a stationary solution to the phase equation;
if instead $\lambda_1=0$, 
then every element $\yo$ of the \omegalimit\ satisfies a problem 
containing a real function~$\mui$ related to the chemical potential~$\mu$.
Such a function $\mui$ is nonunique and time dependent, in general,
as we show by an example.
However, we give sufficient conditions for $\mui$ to be uniquely determined and constant.

\vskip3mm
\noindent {\bf Key words:}
Fractional operators, Cahn--Hilliard systems, longtime \bhv. 
\vskip3mm
\noindent {\bf AMS (MOS) Subject Classification:} 35K45, 35K90, 35R11, \pier{35B40}.
\End{abstract}
\salta
\pagestyle{myheadings}
\newcommand\testopari{\sc Colli \ --- \ Gilardi \ --- \ Sprekels}
\newcommand\testodispari{\sc Generalized fractional Cahn--Hilliard system}
\markboth{\testopari}{\testodispari}
\finqui
%
%%%%%%%%%%%%%%%%%%%%%%%%%%%%%%%%%
%% very beginning
%%%%%%%%%%%%%%%%%%%%%%%%%%%%%%%%%
%%%%%%%%%%%%%%%%%%%%%%%%%%%%%%%%%
%% you may adjust the baseline
%%%%%%%%%%%%%%%%%%%%%%%%%%%%%%%%%

\section{Introduction}
\label{Intro}
\setcounter{equation}{0}

The paper \cite{CGS18} investigates the abstract evolutionary system
\Bsist
  && \dt y + A^{2r} \mu = 0,
  \label{Iprima}
  \\
  && \tau \dt y + B^{2\sigma} y + f'(y) = \mu + u,
  \label{Iseconda}
  \\
  && y(0) = \yz,
  \label{Icauchy}
\Esist
\Accorpa\Ipbl Iprima Icauchy
where $A^{2r}$ and $B^{2\sigma}$, with $r>0$ and $\sigma>0$,
denote fractional powers in the spectral sense of the unbounded linear operators $A$ and~$B$, respectively,
which are supposed to be densely defined in $H:=\Ldue$, with $\Omega\subset\erre^3$,
selfadjoint, and monotone. 
The above system is a generalization
of the Cahn--Hilliard system 
(namely, the \emph{nonviscous} system or the \emph{viscous} one, \juerg{depending on whether} $\tau=0$ or $\tau>0$), 
which models a phase separation process taking place in the container~$\Omega$.  
The unknown functions $\,y\,$ and $\,\mu\,$ stand for the \emph{order parameter}
and the \emph{chemical potential}, respectively,
while $u$ is a given source term. 
Moreover, $\,f\,$ denotes a double-well potential,
for which typical and physically significant examples
are the so-called {\em classical regular potential}, the {\em logarithmic double-well potential\/},
and the {\em double obstacle potential\/}, which are given, in this order,~by
\Bsist
  && f_{reg}(r) := \frac 14 \, (r^2-1)^2 \,,
  \quad r \in \erre, 
  \label{regpot}
  \\
  && f_{log}(r) := \bigl( (1+r)\ln (1+r)+(1-r)\ln (1-r) \bigr) - c_1 r^2 \,,
  \quad r \in (-1,1),
  \label{logpot}
  \\[1mm]
  && f_{2obs}(r) := - c_2 r^2 
  \quad \hbox{if $|r|\leq1$}
  \aand
  f_{2obs}(r) := +\infty
  \quad \hbox{if $|r|>1$}.
  \label{obspot}
\Esist
Here, the constants $c_i$ in \eqref{logpot} and \eqref{obspot} satisfy
$c_1>1$ and $c_2>0$, so that $f_{log}$ and $f_{2obs}$ are nonconvex.
In cases like \eqref{obspot}, one has to split $f$ into a nondifferentiable convex part~$\Beta$ 
(the~indicator function of $[-1,1]$, in the present example) and a smooth perturbation~$\Pi$.
Accordingly, one has to replace the derivative of the convex part
by the subdifferential and interpret \eqref{Iseconda} as a differential inclusion
or, equivalently, as a variational inequality involving $\Beta$ rather than its subdifferential. 
Actually, the latter has been done in~\cite{CGS18},
and \juerg{we do the same} in this paper.

Fractional versions of the Cahn--Hilliard system have been considered by several authors
and are the subject of a number of recent papers.
As for references regarding well-posedness and related problems, 
a rather large list of citations is given in~\cite{CGS18}.
Here we recall some literature dealing with the asymptotic behavior of the solutions.
Indeed, one can find a number of results in this direction
both for the standard Cahn--Hilliard equations and for variants of them.
The latter are obtained, e.g., 
by adding viscosity or memory contributions as well as convective terms;
another possibility is coupling \Ipbl\ with other equations, 
like heat type equations or fluid dynamics equations,
or introducing non--local--in--space terms;
finally, one can replace the classical Neumann boundary conditions by other ones,
e.g., the dynamic boundary conditions.
Without any claim of completeness, 
by starting from~\cite{Zheng},
\pier{we can quote, e.g.,
\cite{AbWi, AkSS1, ChFaPr, CGLN, CGPS2, CGS17, GiMiSchi, GiRo, GS, GrPeSch, JiWuZh, PrVeZa, WaWu, WuZh} 
for the study of the trajectories and related topics,
and 
\cite{EfGaZe, EfMiZe, Gal1, Gal3, Gal2, GaGr1, GaGr2, GrPeSch, LiZho, Mir3, Mir2, Mir1, MiZe2, MiZe1, Seg, ZhaoLiu}}
for the existence of global or exponential attractors and their properties.
However, if \juerg{nonlocal} terms are considered in these papers,
they are not defined as fractional powers in the spectral sense of the operators involved.
On the contrary, our framework is followed in~\cite{CG},
where the longtime \bhv\ of the solutions to a fractional version of the Allen--Cahn \pier{equation} is studied.

Let us come to the content of this paper.
Our aim is studying the \omegalimit\ (in~a suitable topology) of the component $y$ of the solution
to a proper weak version of problem \Ipbl.
The characterization we give (Theorem~\ref{Longtime}) depends 
on the first eigenvalue $\lambda_1$ of the operator~$A$.
If $\lambda_1>0$, then $\mu(t)$ tends to zero as $t$ \juerg{approaches} infinity,
and every element $\yo$ of the \omegalimit\
is a stationary solution, i.e., it solves the equation
\Beq
  B^{2\sigma} \yo + f'(\yo) = \ui,
  \label{Iyostat}
\Eeq
at least in a weak sense,
where $\ui$ is the limit of $u(t)$ as $t$ tends to infinity.
If, instead, $\lambda_1=0$, \juerg{then} the element $\yo$ satisfies a weaker property, namely,
a~weak form of the equation
\Beq
  B^{2\sigma} \yo + f'(\yo) = \ui + \mui(t)
  \quad \Aat,
  \label{Icharomega}
\Eeq
for at least one function $\mui\in L^\infty_{loc}([0,+\infty))$.
We also show that, in the general case, the characterization \eqref{Icharomega} is the best possible
(see Example~\ref{Best}): 
$\mui$~is \juerg{nonconstant and nonunique}, in general, 
and $\mu(t)$ does not converge at infinity.
\pier{On the other hand}, we give sufficient conditions on $f$ and on the solution that ensure that
the function $\mui$ is unique and constant 
and that \eqref{Icharomega} holds in the strong sense (see Proposition~\ref{Goodmui}).

%%%%%%%%%%%%%%%%%%%%%%%%%%%%%%%%%%%%%%%%%%%%%%%%%%%%%%%%%%%%%%%%%%%%%%%%

\section{Statement of the problem and results}
\label{STATEMENT}
\setcounter{equation}{0}

In this section, we state precise assumptions and notations and present our results.
Our framework is the same as in~\cite{CGS18},
and we briefly recall it here, for the reader's convenience.
First of all, the open set $\Omega\subset\erre^3$ is assumed to be bounded, connected and smooth.
We use the notation
\Beq
  H := \Ldue
  \label{defH}
\Eeq
and denote by $\norma\cpto$ and $(\cpto,\cpto)$ the standard norm and inner product of~$H$.
As for the operators involved in our system, we postulate that
\begin{align}
  & A:D(A)\subset H\to H
  \aand
  B:D(B)\subset H\to H
  \quad \hbox{are}
  \nonumber
  \\
  & \hbox{unbounded, monotone, selfadjoint linear operators with compact resolvents.} 
  \label{hpAB} 
\end{align}
We denote by 
$\{\lambda_j\}$ and $\{\lambda'_j\}$ the nondecreasing sequences of the eigenvalues
and by $\{e_j\}$ and $\{e'_j\}$ the (complete) systems of the corresponding orthonormal eigenvectors,
that~is,
\begin{align}
  & A e_j = \lambda_j e_j, \quad
  B e'_j = \lambda'_j e'_j,
  \aand
  (e_i,e_j) = (e'_i,e'_j) = \delta_{ij},
  \quad \hbox{for $i,j=1,2,\dots$},
%  \qquad
  \label{eigen}
  \\
  & 0 \leq \lambda_1 \leq \lambda_2 \leq \dots
  \aand
  0 \leq \lambda'_1 \leq \lambda'_2 \leq \dots,
  \label{eigenvalues}
%  \\
%  && 
  \quad \hbox{with } \,
  \lim_{j\to\infty} \lambda_j
  = \lim_{j\to\infty} \lambda'_j
  = + \infty .
\end{align}
The power $A^r$ of $A$ with an arbitrary positive real exponent~$r$ is given~by
\Beq
  A^r v = \somma j1\infty \lambda_j^r (v,e_j) e_j
  \quad \hbox{for $v\in\VA r$},
  \label{defAr}
\Eeq  
where  
\Beq
  \VA r := D(A^r)
  = \Bigl\{ v\in H:\ \somma j1\infty |\lambda_j^r (v,e_j)|^2 < +\infty \Bigr\}.
  \label{defdomAr}
\Eeq
In principle, we could endow $\VA r$ with the standard graph norm
in order to make $\VA r$ a Hilbert space.
However, we will choose an equivalent Hilbert structure later on.
In the same way, for $\sigma>0$,
we define the power $B^\sigma$ of~$B$.
For its domain, we use the notation
\begin{align}
  & \VB\sigma := D(B^\sigma),
  \quad \hbox{with the norm $\norma\cpto_{B,\sigma}$ associated to the inner product}
%  \label{defBs}
  \non
  \\
  & (v,w)_{B,\sigma} := (v,w) + (B^\sigma v,B^\sigma w)
  \quad \hbox{for $v,w\in \VB\sigma$}.
  \label{defprodBs}
\end{align}
Accordingly, we introduce a space with a negative exponent.
We set
\Beq
  \VA{-r} := (\VA r)^* 
  \quad \hbox{for $r>0$}
  \label{defVAneg}
\Eeq
and use the symbol $\<\cpto,\cpto>_{A,r}$ for the duality pairing
between $\VA{-r}$ and~$\VA r$.
We also identify $H$ with a subspace of $\VA{-r}$
in the usual way, i.e., \juerg{such} that
\Beq
  \< v,w >_{A,r} = (v,w)
  \quad \hbox{for every $v\in H$ and $w\in\VA r$}.
  \label{identification}
\Eeq

At this point, we can start listing our assumptions.
First of all,
\Beq
  \hbox{$r$ and $\sigma$ are fixed positive real numbers and $\tau\in[0,1]$ is fixed as well.}
  \label{hprst}
\Eeq
As for the linear operators, we postulate, besides~\eqref{hpAB}, that
\Bsist
  && \hbox{either} \quad
  \lambda_1 > 0 
  \quad \hbox{or} \quad
  \hbox{$0=\lambda_1<\lambda_2$ and $e_1$ is a constant;}
  \label{hpsimple}
  \\
  && \hbox{if\quad $\lambda_1=0$, \juerg{then}
  the constant functions belong to $\VB\sigma$}.
  \label{hpVB}
\Esist
In \cite{CGS18} some remarks are given on the above assumptions.
Moreover, it is shown that an equivalent Hilbert structure on $\VA r$ is obtained
by taking the norm defined~by
\Beq
  \norma v_{A,r}^2 := \left\{ 
  \begin{aligned}
  & \norma{A^r v}^2
  = \somma j1\infty |\lambda_j^r (v,e_j)|^2
  \qquad \hbox{if $\lambda_1>0$,}
  \\
  & |(v,e_1)|^2 + \norma{A^r v}^2
  = |(v,e_1)|^2 + \somma j2\infty |\lambda_j^r (v,e_j)|^2
  \qquad \hbox{if $\lambda_1=0$}.
  \end{aligned}
  \right.
  \label{defnormaAr}
\Eeq
We notice that the term $(v,e_1)$ appearing in \eqref{defnormaAr} in the case $\lambda_1=0$
is proportional to the mean value of~$v$
\Beq
  \mean v := \frac 1 {|\Omega|} \iO v \,,
  \label{defmean}
\Eeq
since $e_1$ is a constant by~\eqref{hpsimple}.
In particular, we have the Poincar\'e type inequality
\Beq
  \norma v \leq C_P \, \norma{A^r v}
  \quad \hbox{for every $v\in\VA r$ with $\mean v=0$, \ if $\lambda_1=0$}.
  \label{poincare}
\Eeq
For the nonlinearity $f$ appearing in our system, we split it as $f=\Beta+\Pi$ 
and postulate the following properties
(which are fulfilled by all of the important potentials \accorpa{regpot}{obspot}):
\Bsist
  & \Beta : \erre \to [0,+\infty]
  \quad \hbox{is convex, proper, and l.s.c., \ with} \quad
  \Beta(0) = 0;
  \label{hpBeta}
  \\
  \separa
  & \Pi : \erre \to \erre
  \quad \hbox{is of class $C^1$ with a \Lip\ continuous first derivative;}
  \label{hpPi}
  \\
  & \mbox{it holds }\,\,\displaystyle \liminf_{|s|\nearrow+\infty} \displaystyle \frac {\Beta(s)+\Pi(s)} {s^2} > 0.
  \label{hpcoerc}
\Esist
We set, for convenience,
\Beq
  \beta := \partial\Beta , \quad
  \pi := \Pi' , \quad
  \Lpi = \hbox{the \Lip\ constant of $\pi$,}
  \aand
  \Lpi' := \Lpi + 1 \,.
  \label{defbetapi}  
\Eeq
Moreover, we term $D(\Beta)$ and $D(\beta)$ the effective domains of $\Beta$ and~$\beta$, respectively,
and notice that $\beta$ is a maximal monotone graph in $\erre\times\erre$.

At this point, we can state the problem under investigation,
and we do it on the half-line $t\geq0$,
due to the subject of the present paper.
The data are required to satisfy
\Bsist
  && u \in W^{1,1}_{loc}([0,+\infty);H)
  \aand 
  \dt u\in\LL1H.
  \label{hpu}
  \\
  && \yz \in \VB\sigma
  \aand
  \Beta(\yz) \in \Luno . 
  \label{hpyz}
  \\
  && \hbox{If $\lambda_1=0$, \juerg{then}} \,\,\,\,
  \mz := \mean\yz
  \,\,\,\, \hbox{belongs to the interior of $D(\beta)$} .
  \label{hpmz}
\Esist
\Accorpa\HPdati hpu hpmz
A~solution to our system is a pair $(y,\mu)$ fulfilling the regularity requirements
\Bsist
  && y \in \L\infty{\VB\sigma}, \quad
  \dt y \in \L2{\VA{-r}}
  \aand
  \tau \dt y \in \L2H,
  \qquad
  \label{regy}
  \\
  && \mu \in \L2{\VA r},
  \label{regmu}
  \\
  && \Beta(y) \in L^1(\Omega\times(0,T)),
  \label{regBetay}
\Esist
\Accorpa\Regsoluz regy regBetay
for every $T>0$,
and satisfying the following weak formulation of the equations \accorpa{Iprima}{Icauchy}:
\Bsist
  && \< \dt y(t) , v >_{A,r}
  + ( A^r \mu(t) , A^r v )
  = 0
  \quad \hbox{for every $v\in\VA r$\, and \Aat},
  \qquad\quad
  \label{prima}
  \\[2mm]
  && \juerg {\bigl(\tau } \dt y(t) , y(t) - v \bigr)
  + \bigl( B^\sigma y(t) , B^\sigma( y(t)-v) \bigr)
  \non
  \\
  && \quad {}
  + \iO \Beta(y(t))
  + \bigl( \pi(y(t)) - u(t) ,  y(t)-v \bigr)
  \leq \bigl( \mu(t) ,  y(t)-v \bigr)
  + \iO \Beta(v)
  \non
  \\
  && \quad \hbox{for every $v\in\VB\sigma$\, and \Aat},
  \label{seconda}
  \\
  && y(0) = \yz \,.
  \label{cauchy}
\Esist
\Accorpa\Pbl prima cauchy
We remark that, if $\lambda_1=0$, then $A^r(1)=0$ by \eqref{hpsimple},
so that \eqref{prima} implies that
\Bsist
  && \frac d{dt} \iO y(t) = 0
  \quad \Aat, \quad
  \hbox{i.e.,} \quad
  \non
  \\
  && \mean y(t) = \mz
  \quad \hbox{for every $t\in[0,+\infty)$}.
  \label{conservation}
\Esist
The well-posedness result stated below was proved in \cite{CGS18} under a different assumption on~$u$.
Namely, in studying the problem on the finite time interval~$(0,T)$, it was assumed that $u\in\H1H$,
while \eqref{hpu} only implies that $u\in\W{1,1}H$.
However, we point out that our assumption is sufficient to obtain the same result.
We will give some explanation on this in the next section.

\Bthm
\label{Wellposedness}
Let the assumptions \eqref{hpAB}, \accorpa{hprst}{hpVB} and 
\accorpa{hpBeta}{hpcoerc} 
on the structure of the system,
and \HPdati\ on the data, be fulfilled.
Then there exists a pair $(y,\mu)$ satisfying \Regsoluz\
and solving problem \Pbl.
Moreover, the component $y$ of the solution is \juerg{uniquely determined}.
\Ethm

In \cite[Rem.~4.1]{CGS18}, sufficient conditions \juerg{were given that ensure uniqueness also for $\mu$.}
However, the aim of this paper is the study of the longtime \bhv\ of the \juerg{component~$y$ alone}.
The rather weak regularity conditions \eqref{regy} imply that
\Beq
  y :[0,+\infty) \to \VB\sigma
  \quad \hbox{is weakly continuous}.
  \non
\Eeq
This enables us to the define the following (possibly empty) \omegalimit\ \pier{set}
\Beq
  \omega = \omega (\yz,u) 
  := \graffe{\yo\in\VB\sigma:\ y(\tn)\to \yo \quad \hbox{weakly in $\VB\sigma$ for some $\graffe{\tn}\nearrow+\infty$}
  }.
  \label{defomega}
\Eeq

Here is our result, which holds under the additional assumption that 
$u(t)$ has a limit $\ui$ as $t$ tends to infinity 
in the sense of the forthcoming~\eqref{hpui}.
The second part of the statement distinguishes two cases regarding the first eigenvalue $\lambda_1$ of~$A$.
If $\lambda_1$ is positive, \juerg{then} every element of the \omegalimit\ is a stationary solution
in the sense specified below;
if instead $\lambda_1=0$, \juerg{then} the elements of the \omegalimit\ just satisfy a weaker property.

\Bthm
\label{Longtime}
Let the assumptions \eqref{hpAB}, \accorpa{hprst}{hpVB} and 
\accorpa{hpBeta}{hpcoerc} 
on the structure of the system,
and \HPdati\ on the data, be fulfilled.
In addition, assume that there \juerg{is some} $\ui\in H$ such that
\Beq
  u - \ui \in \LL2H,
  \label{hpui}
\Eeq
and let $(y,\mu)$ be a solution to \Pbl\ according to Theorem~\ref{Wellposedness}.
Then the \omegalimit\ \eqref{defomega} is \juerg{nonempty}.
Moreover, it is characterized as follows:

\noindent
$i)$ If $\lambda_1>0$, \juerg{then} every element $\yo\in\omega$ satisfies
\Bsist
  && \bigl( B^\sigma \yo , B^\sigma (\yo-v) \bigr) 
  + \iO \Beta(\yo) 
  + \bigl( \pi(\yo) - \ui , \yo-v \bigr)
  \leq \iO \Beta(v)
  \non
  \\
  && \quad \hbox{for every $v\in\VB\sigma$}.
  \label{yostat}
\Esist
$ii)$ If $\lambda_1=0$,
then, for every element $\yo\in\omega$, 
there exists \juerg{some} $\mui\in L^\infty_{loc}([0,+\infty))$ such~that
\Bsist
  && \bigl( B^\sigma \yo , B^\sigma (\yo-v) \bigr) 
  + \iO \Beta(\yo) 
  + \bigl( \pi(\yo) - \ui , \yo-v \bigr)
  \non
  \\
  && \leq \bigl( \mui(t) , \yo-v \bigr)
  + \iO \Beta(v)
  \non
  \\
  && \quad \hbox{for every $v\in\VB\sigma$ and \Aat}.
  \label{charomega}
\Esist
\Ethm

In \eqref{hpui}, $\ui$ obviously denotes the function $[0,+\infty)\ni t\mapsto\ui\in H$
rather than the element $\ui\in H$.
In the \rhs\ of~\eqref{charomega}, $\mui(t)$ denotes the constant function 
$\Omega\ni x\mapsto\mui(t)$ 
rather than the real value~$\mui(t)$.
Conventions of this type \juerg{will be} used also in the following.

The part $ii)$ of the above result seems \juerg{to be} rather poor. 
Nevertheless, \juerg{this} characterization is the best possible for the general case, that is,
one can \juerg{neither} expect uniqueness for~$\mui$, 
nor further properties for it, as the following example shows.
Notice that assuming that $A$ and $B$ are particularly good operators does not help at all.

\Bex
\label{Best}
Let the operators $A$ and $B$ satisfy the hypotheses of Theorem~\ref{Longtime},
and assume that $\lambda_1=0$.
Moreover, let us choose $\Pi=0$ and $\Beta$ given~by
\Beq
  \Beta(r) := r^2 + |r|
  \quad \hbox{for $r\in\erre$}.
  \non
\Eeq
Then \accorpa{hpBeta}{hpcoerc} are satisfied.
But $\beta$ is multivalued, since $\beta(0)=\sign(0)=[-1,1]$.
Thus, if we take $\yz=0$, $u=0$, and any function $\bar\mu\in L^\infty(0,+\infty)$
satisfying $|\bar\mu(t)|\leq1$ \Aat, \juerg{then}
a~solution $(y,\mu)$ to problem \Pbl\ is given by the formulas
$y(x,t)=0$ and $\mu(x,t)=\bar\mu(t)$.
Indeed, $(y,\mu)$ trivially solves the first equation~\eqref{prima} 
(since $\mu$ is space independent), as~well as~\eqref{cauchy};
moreover, the variational inequality \eqref{seconda} is solved in the stronger form
\Beq
  \tau \dt y + B^{2\sigma} y + \xi + \pi(y) = \mu + u
  \quad \hbox{with} \quad
  \xi \in \beta(y),
  \non
\Eeq
since we can take $\xi=\mu$
(we have $\bar\mu(t)\in[-1,1]=\beta(0)$, indeed).
So, the only element $\yo$ of the \omegalimit\ is $\yo=0$,
while we \pier{have lots of possible $\mui$'s, namely, the set of such functions 
\juerg{coincides} with the set of the admissible functions termed $\bar\mu $ before.}
\Eex

On the contrary, 
under further conditions on $\beta$ and on the solution,
the characterization in the case $\lambda_1=0$ can be improved.
Here are the new requirements:
\Bsist
  && \hbox{$D(\beta)$ is an open interval, and $\beta$ is a single-valued $C^1$ function.}
  \label{regbeta}
  \\[2mm]
  && \hbox{There exists a compact interval $[a,b]\subset D(\beta)$ such that}
  \nonumber
  \\
  && \hbox{\quad $y(x,t)\in[a,b]$ for a.a.\ $(x,t)\in\Omega\times(0,+\infty)$.}
  \label{gb}
  \\[2mm]
  && \hbox{$\VB\sigma\cap\Linfty$ \ is dense in \ $\VB\sigma$}.
  \label{density}
\Esist
The above assumptions (\pier{with} \eqref{gb} only in a given finite time interval $(0,T)$)
\juerg{have been} introduced in the paper~\cite{CGS19}.
One of the motivations \juerg{was} the derivation of the strong form of \eqref{seconda},~i.e.,
\Beq
  \tau \, \dt y + B^{2\sigma} y + \beta(y) + \pi(y) = \mu + u \,.
  \label{strong}
\Eeq
Precisely, it \juerg{has been} proved that $y\in\L2{\VB{2\sigma}}$ and that \eqref{strong} is satisfied almost everywhere
(see Rem.~3.5 and the subsequent lines of~\cite{CGS19}, where some comments on \accorpa{regbeta}{density} 
were given as well).
Here, we point out that the proof of the derivation of \eqref{strong} 
also holds true for the half--line $t\geq0$ if \eqref{gb} is assumed.
We use \accorpa{regbeta}{density} in the result stated below.

\Bprop
\label{Goodmui}
In addition to the assumptions of Theorem~\ref{Longtime},
suppose that \accorpa{regbeta}{density} are satisfied.
Then the function $\mui$ appearing in \eqref{charomega} is uniquely determined  and constant.
Moreover, $\yo\in\VB{2\sigma}$, and the pair $(\yo,\mui)$ satisfies the equation
\Beq
  B^{2\sigma} \yo + \beta(\yo) + \pi(\yo) = \mui + \ui 
  \quad \aeO .
  \label{strongi}
\Eeq
\Eprop

Theorem~\ref{Longtime} and Proposition~\ref{Goodmui} \juerg{will be} proved in the last section.
In the next one, we establish some auxiliary global estimates.
To this \juerg{end}, we also recall the approximation and the discretization of problem \Pbl\
given in~\cite{CGS18}.

\Bnot
\label{Notation}
In the \juerg{remainder} of the paper, we will use the same small letter $c$
for (possibly) different constants that depend only 
on the structure of our system (but~$\tau$) and on the assumptions on the data.
When some final time $T$ is considered,
the symbol $c_T$ denotes (possibly different) constants
that depend on~$T$ in addition.
On the contrary, precise constants we could refer to are treated in a different way
(see, e.g., the forthcoming~\eqref{coerclambda}, where greek and capital letters are used).
\Enot

%%%%%%%%%%%%%%%%%%%%%%%%%%%%%%%%%%%%%%%%%%%%%%%%%%%%%%%%%%%%%%%%%%%%%%%%

\section{Global estimates}
\label{GLOBAL}
\setcounter{equation}{0}

The proof of Theorem~\ref{Longtime} is based on some global--in--time a~priori estimates,
which we derive in this section by starting from the approximating and discrete problems introduced in~\cite{CGS18}.
Thus, some recalls are needed.

The approximation of problem \Pbl\ by a more regular one
relies on the use of the Moreau--Yosida regularizations 
$\Betal$ and $\betal$ \,of\, $\Beta$ \juerg{and} $\beta$ at the level $\lambda>0$
(see, e.g., \cite[p.~28 and p.~39]{Brezis}).
We notice that, by~accounting for \eqref{hpcoerc}, the inequalities
\Beq
  \Betal(s) + \Pi(s)
  \geq \alpha \, s^2 - C
  \geq - C'
  \label{coerclambda}
\Eeq
hold true for some \juerg{positive} constants $\alpha,C, C'$, every $s\in\erre$, and every 
\juerg{sufficiently small $\lambda>0$}. \pier{In case the reader aims to check  \eqref{coerclambda}, we suggest the use of the following representation of $\Betal$, namely
\begin{gather}
	\widehat{\beta }_{\lambda }(s)
	:=\inf_{r \in \mathbb{R}}\left\{ \frac{1}{2\lambda } |r-s|^2
	+\widehat{\beta }(r) \right\} 
	= 
	\frac{1}{2\lambda } 
	\bigl| s-J_\lambda  (s) \bigr|^2+\widehat{\beta }\bigl (J_\lambda (s) \bigr ),
	\label{pier1}
\end{gather}
where 
$J_\lambda :\mathbb{R} \to \mathbb{R}$  denotes the resolvent operator associated to $\beta$, that is, 
$ J_\lambda (s)$ is defined as the unique solution \juerg{to} the multi-equation
$$J_\lambda (s) + \lambda \beta ( J_\lambda (s)) \ni s \equiv J_\lambda (s) + \lambda \betal( s) \quad \hbox{ for all } \, s\in \mathbb{R}.$$ 
Indeed, by combining \eqref{hpBeta}--\eqref{hpcoerc}, which imply
\Beq
  \Beta(s) + \Pi(s)
  \geq 2 \alpha \, s^2 - c \quad \hbox{ for all } \, s\in \mathbb{R}
   \non
\Eeq
and for some constant $\alpha >0$, along with \eqref{pier1} and the Taylor formula with integral remainder to estimate the difference $\Pi (s) - \Pi ( J_\lambda (s)) $, one can arrive at 
\eqref{coerclambda}.}

The approximating problem on any finite time integral $(0,T)$ is obtained
by replacing $\Beta$ in \eqref{seconda} by~$\Betal$,
namely,
\Bsist
  && \< \dt\yl(t) , v >_{A,r}
  + ( A^r \mul(t) , A^r v )
  = 0
  \quad \hbox{for every $v\in\VA r$ and \aat},
  \qquad
  \label{primal}
  \\[1mm]
  && \juerg{\bigl(\tau} \dt\yl(t) , \yl(t) - v \bigr) 
  + \bigl( B^\sigma\yl(t) , B^\sigma(\yl(t)-v) \bigr)
  \non
  \\
  && \quad {}
  + \iO \Betal(\yl(t))
  + \bigl( \pi(\yl(t)) - u(t) ,  \yl(t)-v \bigr)
  \non
  \\
  && \leq \bigl( \mul(t) ,  \yl(t)-v \bigr)
  + \iO \Betal(v)
  \quad \hbox{for every $v\in\VB\sigma$ and \aat},
  \label{secondal}
  \\[1mm]
  && \yl(0) = \yz .
  \label{cauchyl}
\Esist
\Accorpa\Pbll primal cauchyl
In principle, the regularity required for the solution $(\yl,\mul)$ is still \juerg{given by} \Regsoluz.
However, due to the \Lip\ continuity of~$\betal$, \eqref{regBetay}~can be improved.
Namely, \eqref{regy}~implies $\betal(\yl)\in\L2H$. 
Using this and the fact that $\Betal$ is differentiable and $\betal$ is its derivative,
one sees that, in place of \eqref{secondal}, 
one can equivalently consider the pointwise variational equation 
\begin{align}
  & \pier{\bigl(\tau \dt\yl(t) ,  v \bigr)+{}}
  \bigl( B^\sigma\yl(t) , B^\sigma v \bigr)
  + \bigl( \betal(\yl(t)) + \pi(\yl(t)) - u(t) ,  v \bigr)
  = \bigl( \mul(t) ,  v \bigr)
  \non
  \\[1mm]
  & \quad \hbox{for every $v\in\VB\sigma$ and \aat}.
  \label{eqsecondal}
\end{align}
In \cite{CGS18}, it was shown that the above problem is~well-posed
and~that its unique solution $(\yl,\mul)$ converges
to a solution $(y,\mu)$ to problem \Pbl\
in the weak topology associated with the regularity requirements, essentially.
Moreover, the solution $(\yl,\mul)$ is obtained as the limit of suitable interpolant functions
constructed by starting from the solution to a proper discrete problem.
For the reader's convenience, we recall both the notation for the interpolants and the discrete problem.

Let $N$ be a positive integer and $Z$ be one of the spaces $H$, $\VA r$,~$\VB\sigma$.
We set $h:=T/N$ and $I_n:=((n-1)h,nh)$ for $n=1,\dots,N$.
Given $z=(z_0,z_1,\dots ,z_N)\in Z^{N+1}$,
the piecewise constant and piecewise linear interpolants 
\Beq
  \overline z_h \in \L\infty Z , \quad
  \underline z_h \in \L\infty Z 
  \aand
  \hat z_h \in \W{1,\infty}Z
  \non
\Eeq
are defined by setting 
\Bsist
  && \hskip -2em
  \overline z_h(t) = z^n
  \aand
  \underline z_h(t) = z^{n-1}
  \quad \hbox{for a.a.\ $t\in I_n$, \ $n=1,\dots,N$},
  \label{pwconstant}
  \\
  && \hskip -2em
  \hat z_h(0) = z_0
  \aand
  \dt\hat z_h(t) = \frac {z^{n+1}-z^n} h
  \quad \hbox{for a.a.\ $t\in I_n$, \ $n=1,\dots,N$}.
  \qquad
  \label{pwlinear}
\Esist
The discrete problem consists in finding
two $(N+1)$-tuples $(y^0,\dots,y^N)$ and $(\mu^0,\dots,\mu^N)$
satisfying
\Bsist
  y^0 = \yz \,, \quad
  \mu^0 = 0 , \quad
  (\yu,\dots,y^N) \in (\VB{2\sigma})^N
  \aand
  (\muu,\dots,\mu^N) \in (\VA{2r})^N,
  \label{regdiscr}
\Esist
and solving
\begin{align}
  & \dhyn + \munp + A^{2r} \munp
  \,=\, \mun,
  \label{primad}
  \\[1mm]
  & \tau \, \dhyn
  + (\Lpi' I + B^{2\sigma} + \betal + \pi)(\ynp)
   \,=\, \Lpi' \yn + \munp + \unp,
  \label{secondad}
\end{align}
for $n=0,1,\dots,N-1$, where $I:H\to H$ is the identity, $\Lpi'$~is given by 
\eqref{defbetapi}, and
\Beq
  \un := u(nh)
  \quad \hbox{for $n=0,1,\dots,N$}.
  \label{defun}
\Eeq
Precisely, it \juerg{has been} proved that such a discrete problem is uniquely solvable.
Moreover, as just said, some of the interpolants defined above by starting
from the discete solution converge to the solution $(\yl,\mul)$ to the regularized problem \Pbll.

Now, we start estimating.
It is understood that the assumptions of Theorem~\ref{Longtime} are in force.
In particular, every constant $c$ we introduce
will depend only on \juerg{these} assumptions.
We closely follow \juerg{the lines of} \cite{CGS18}.
However, we modify the argument a little and obtain estimates that are uniform with respect to~$T$.
In doing this modification, we also avoid using the regularity condition $\dt u\in\L2H$\juerg{, which was supposed in~\cite{CGS18},}
and just owe to the regularity $\dt u\in\L1H$ 
(but uniformly with respect to $T$ in the sense of \eqref{hpu}
in order to obtain a global--in--time estimate).
Since this is the only point of \cite{CGS18} where the $\L2H$ regularity for~$\dt u$ is accounted for,
the well-posedness \pier{result in} Theorem~\ref{Wellposedness} holds under our assumption~\eqref{hpu}, 
as announced before its statement.

\step
First uniform estimate

We test \eqref{primad} and~\eqref{secondad} 
(by~taking the scalar product in~$H$)
by~$h\munp$ and $\ynp-\yn$, respectively,
and add the resulting identities.
Noting an obvious cancellation, we obtain the equation
\begin{align}
  & h (\munp-\mun,\munp)
  + h (A^{2r} \munp,\munp)
    + \frac \tau h \, \norma{\ynp-\yn}^2
	\non\\[1mm]	
  &+ (B^{2\sigma}\ynp,\ynp-\yn)
    + \bigl( (\Lpi' I + \betal + \pi) (\ynp),\ynp-\yn \bigr)
  \non
  \\[1mm]
  & =\, \Lpi' (\yn,\ynp-\yn)
  + (\unp,\ynp-\yn).
  \non
\end{align}
Now, we observe that the function 
$\,r\mapsto\frac{\Lpi'}2\,r^2+\Betal(r)+\Pi(r)\,$
is convex on~$\erre$, since $\Betal$ is convex and $\,|\pi'|\leq\Lpi$.
Thus, we have that
\Bsist
  && \bigl( (\Lpi' I + \betal + \pi)(\ynp) , \ynp-\yn \bigr)
  \non
  \\  
  && \geq \frac {\Lpi'} 2 \, \norma\ynp^2 + \iO \bigl( \Betal(\ynp) + \Pi(\ynp) \bigr)
  - \frac {\Lpi'} 2 \, \norma\yn^2 - \iO \bigl( \Betal(\yn) + \Pi(\yn) \bigr).
  \non
\Esist
We easily deduce~that
\Bsist
  && \frac h2 \, \norma\munp^2
  + \frac h2 \, \norma{\munp-\mun}^2
  - \frac h2 \, \norma\mun^2
  + h \norma{A^r\munp}^2
  \non
  \\
  && \quad {}
  + \frac \tau h \, \norma{\ynp-\yn}^2
  + \frac 12 \, \norma{B^\sigma \ynp}^2
  + \frac 12 \, \norma{B^\sigma(\ynp-\yn)}^2
  - \frac 12 \, \norma{B^\sigma\yn}^2
  \non
  \\
  && \quad {}
  + \frac {\Lpi'} 2 \, \norma\ynp^2 + \iO \bigl( \Betal(\ynp) + \Pi(\ynp) \bigr)
  - \frac {\Lpi'} 2 \, \norma\yn^2 - \iO \bigl( \Betal(\yn) + \Pi(\yn) \bigr)
  \non
  \\
  \separa
  && \leq - \frac {\Lpi'} 2 \bigl(
       \norma\yn^2
       - \norma\ynp^2
       + \norma{\ynp-\yn}^2
  \bigr)
  + (\unp , \ynp-\yn).
  \non
\Esist
Then, we first rearrange and then sum up for $n=0,\dots,k-1$ with $k\leq N$,
employing summation by parts in the last term.
We thus arrive at the inequality
\Bsist
  && \frac h2 \, \norma\muk^2
  + \somma n0{k-1} \frac h2 \, \norma{\munp-\mun}^2
  + \somma n0{k-1} h \norma{A^r\munp}^2
  \non
  \\
  && \quad {}
  + \tau \somma n0{k-1} h \Norma{\dhyn}^2
  + \frac 12 \, \norma{B^\sigma \yk}^2
  - \frac 12 \, \norma{B^\sigma \yz}^2
  + \somma n0{k-1} \frac 12 \, \norma{B^\sigma(\ynp-\yn)}^2
  \non
  \\
  && \quad {}
  + \iO \bigl( \Betal(\yk) + \Pi(\yk) \bigr)
  - \iO \bigl( \Beta(\yz) + \Pi(\yz) \bigr)
  + \frac {\Lpi'} 2 \somma n0{k-1} \norma{\ynp-\yn}^2
  \non
  \\
  \separa
  && \leq (\uk,\yk) - (u^1,\yz) - \somma n1{k-1} (\unp-\un,\yn).
  \label{perprimaunif}
\Esist
Next, we observe that \eqref{coerclambda} implies that
\Beq
  \iO \bigl( \Betal(\yk) + \Pi(\yk) \bigr)
  \geq \frac 12 \iO \bigl( \Betal(\yk) + \Pi(\yk) \bigr)
  + \frac \alpha 2 \, \norma\yk^2 - c
  \non
\Eeq
for every sufficiently small $\lambda>0$
and that the integrals are bounded from below.
Moreover, we differently deal with the \rhs\ of~\eqref{perprimaunif} with respect to~\cite{CGS18}.
Namely, we estimate it as follows:
\Bsist
  && (\uk,\yk) - (u^1,\yz) - \somma n1{k-1} (\unp-\un,\yn)
  \non
  \\[-5pt]
  && \leq \frac \alpha 4 \, \norma\yk^2
  + \frac 1\alpha \, \norma\uk^2
  + \norma{u^1} \, \norma\yz
  + \somma n1{k-1} h \, \Norma{\frac{\unp-\un}h} \, \juerg{\|\yn\|} \,.
  \non
\Esist
At this point, we combine \eqref{perprimaunif} with the inequalities just obtained
and apply the discrete Gronwall-Bellman lemma given in \cite[Thm.~1]{Yeh}
by observing that
\Beq
  \norma\uk \leq \norma{u(0)} + \norma{\dt u}_{\LL1H},
  \aand
  \somma n1{k-1} h \, \Norma{\frac{\unp-\un}h}
  \leq \norma{\dt u}_{\LL1H},
  \non
\Eeq
and that the above norm of $\dt u$ is finite by~\eqref{hpu}.
We obtain the estimate
\begin{align}
  & h \, \norma\muk^2
  + \somma n0{k-1} \pier{h} \, \norma{\munp-\mun}^2
  + \somma n0{k-1} h \norma{A^r\munp}^2
  + \tau \somma n0{k-1} h \Norma{\dhyn}^2
  \non
  \\
  & \quad {}
  + \norma\yk_{B,\sigma}^2
  + \somma n0{k-1} \norma{B^\sigma(\ynp-\yn)}^2
  + \iO \bigl( \Betal(\yk) + \Pi(\yk) \bigr)
  + \somma n0{k-1} \norma{\ynp-\yn}^2
  \non
  \\
  & \leq c 
  \quad \hbox{for $k=0,\dots,N$}.
  \label{primastimad}
\end{align}
In terms of the interpolants (see also \cite[Prop.~3.9]{CGS18}), 
by neglecting the first contribution
and recalling that $\mu^0=0$,
we have that
\Bsist
  && \norma{\overmuh-\undermuh}_{\L2H}
  + \norma{A^r \overmuh}_{\L2H}
  + \norma{A^r \undermuh}_{\L2H}
  \non
  \\[2pt] 
  && \quad {}
  + \norma\underyh_{\L\infty{\VB\sigma}}
  + \norma\overyh_{\L\infty{\VB\sigma}}
  + \norma\yh_{\L\infty{\VB\sigma}}
  \non
  \\
  && \quad {}
  + h^{-1/2} \norma{B^\sigma(\overyh-\underyh)}_{\L2H}
  + \tau^{1/2} \norma{\dt\yh}_{\L2H}
  \non
  \\
  && \quad {}
  + \norma{\Betal(\overyh)+\Pi(\overyh)}_{\L\infty\Luno}
  + \pier{h^{-1/2} \norma{\overyh-\yh}_{\L2H}}
  \leq c \,. \qquad
  \label{primaunif}
\Esist

\step
Second uniform estimate

By observing that \eqref{primad} implies 
$\dt\yh+\overmuh+A^{2r}\overmuh=\undermuh$,
whence also
\Bsist
  && \ioT \bigl( \dt\yh(t),v(t) \bigr) \, dt
  = \ioT \bigl( (\undermuh - \overmuh)(t) , v(t) \bigr) \, dt
  - \ioT \bigl( A^r\overmuh(t) ,A^r v(t) \bigr) \, dt
  \non
  \\
  && \leq c \, \bigl(
    \norma{\undermuh - \overmuh}_{\L2H} 
    + \norma{A^r\overmuh}_{\L2H} 
  \bigr) \, \norma v_{\L2{\VA r}}
  \non
\Esist
for every $v\in\L2{\VA r}$,
we deduce that
\Beq
  \norma{\dt\yh}_{\L2{\VA{-r}}}
  \leq c \, \bigl(
    \norma{\undermuh-\overmuh}_{\L2H}
    + \norma{A^r\overmuh}_{\L2H} 
  \bigr) \,.
  \non
\Eeq
Hence, from \eqref{primaunif} we infer that
\Beq
  \norma{\dt\yh}_{\L2{\VA{-r}}} \leq c \,.
  \label{secondaunif}
\Eeq

\step
Basic estimate

We recall that estimates \accorpa{primaunif}{secondaunif}
hold for every $N>1$, every \juerg{sufficiently small $\lambda>0$,}  and every $T>0$.
Now, we owe to the convergence results of~\cite{CGS18}.
We deduce that
\Bsist
  && \norma\yl_{\L\infty{\VB\sigma}}
  + \norma{\dt\yl}_{\L2{\VA{-r}}}
  + \tau^{1/2} \norma{\dt\yl}_{\L2H}
  + \norma{A^r\mul}_{\L2H}
  \leq c \,.
  \non
\Esist
Since $c$ is independent of both $\lambda$ and $T$, \pier{at the limit as $\lambda \searrow 0$} we conclude that
\begin{gather}
  y \in \LL\infty{\VB\sigma}, \quad
  \dt y \in \LL2{\VA{-r}},
%  \non
%  \\
  \quad\hbox{and }\,
  A^r\mu \in \LL2H,
  \label{basic}
  \\[.2cm]
   \dt y \in \LL2H
  \quad \hbox{if $\tau>0$}.
  \label{basictau}
\end{gather}

\section{Longtime \bhv}
\label{LONGTIME}
\setcounter{equation}{0}

This section is devoted to the proofs of our results on the longtime \bhv.
We start with the proof of Theorem~\ref{Longtime}.

\step First part

Since $y$ belongs to $\LL\infty{\VB\sigma}$ by the first \juerg{conclusion of} \eqref{basic},
we deduce that the \omegalimit\ $\omega$ given by \eqref{defomega} is nonempty.
Thus, the first sentence of our result is established. 
Let us come to the second part.

%%%%%%%%%%%%%%%%%%%%%%%%%%%%%%%%

%%%%%%%%%%%%%%%%%%%%%%%%%%%%%%%%

% re-defined macros

\def\yn{y_n}
\def\mun{\mu_n}
\def\un{u_n}

%%%%%%%%%%%%%%%%%%%%%%%%%%%%%%%%

%%%%%%%%%%%%%%%%%%%%%%%%%%%%%%%%

\step Second part, first case

We \juerg{first} assume that $\lambda_1>0$.
We pick an arbitrary element $\yo\in\omega$ and a sequence $\graffe{\tn}$ as in~\eqref{defomega},
and we prove that $\yo$ is a stationary solution in the sense of~\eqref{yostat}. \juerg{To this end, we}
 define the functions $\yn$, $\mun$, and $\un$, on~$(0,+\infty)$ by setting, \Aat,
\Beq
  \yn(t) := y(t+\tn) , \quad
  \mun(t) := \mu(t+\tn),
  \aand
  \un(t) := u(t+\tn).
  \non
\Eeq
We notice that \eqref{hpu} and \eqref{hpui} imply that
\Beq
  \norma\un_{\LL\infty H} \leq c,
  \aand
  \un - \ui \to 0 
  \quad \hbox{strongly in $\LL2H$}.
  \label{stimeun}
\Eeq
Moreover, from \eqref{basic} we clearly deduce that
\Bsist
  && \norma\yn_{\LL\infty{\VB\sigma}} \leq c, 
  \label{ynbdd}
  \\
  && \dt\yn \to 0,
  \quad \hbox{strongly in $\LL2{\VA{-r}}$},
  \label{dtyntozero}
  \\
  && A^r\mun \to 0 
  \quad \hbox{strongly in $\LL2H$,}
  \label{forsemstimeun}  
\Esist
whence also
\Beq
  \mun \to 0
  \quad \hbox{strongly in $\LL2{\VA r}$,}
  \label{stimemun}
\Eeq
since $\lambda_1>0$.
In addition, we have that
\Beq
  \dt\yn \to 0
  \quad \hbox{strongly in $\LL2H$ \quad if $\tau>0$} \,.
  \label{dtyntozerobis}
\Eeq
By weak-star compactness, we deduce from \eqref{ynbdd} that \pier{there exists 
\juerg{some element} $\yi \in 
\LL\infty{\VB\sigma}$ such that}
\Beq
  \yn \to \yi
  \quad \hbox{weakly star in $\LL\infty{\VB\sigma}$},
  \label{convyn}
\Eeq
at least for a (not relabeled) subsequence.
Now, we fix \pier{an arbitrary} time $T>0$ and look for the problem solved by~$\yi$ on~$(0,T)$.
It is clear that $(\yn,\mun)$ satisfies the variational inequality
\Bsist
  && \juerg{\bigl(\tau} \dt\yn(t) , \yn(t) - v \bigr)
  + \bigl( B^\sigma \yn(t) , B^\sigma( \yn(t)-v) \bigr)
  \nonumber
  \\
  && \quad {}
  + \iO \Beta(\yn(t))
  + \bigl( \pi(\yn(t)) - \un(t) ,  \yn(t)-v \bigr)
  \leq \bigl( \mun(t) ,  \yn(t)-v \bigr)
  + \iO \Beta(v)
  \nonumber
  \\
  && \quad \hbox{for every $v\in\VB\sigma$ and \aat,}
  \label{secondan}
\Esist
as well as its integrated version
\Bsist
  &&  \ioT \juerg{\bigl(\tau} \dt\yn(t) , \yn(t) - v(t) \bigr) \, dt 
  + \ioT \bigl( B^\sigma \yn(t) , B^\sigma( \yn(t)-v(t)) \bigr)\,dt 
  \non
  \\ 
  && \quad {}
  + \intQ \Beta(\yn) 
  + \ioT \bigl( \pi(\yn(t)) - \un(t) ,  \yn(t)-v(t) \bigr)\,dt 
  \non
  \\
  && \leq \ioT \bigl( \mun(t) ,  \yn(t)-v(t) \bigr)\,dt
  + \intQ \Beta(v)
  \quad \hbox{for every $v\in\L2{\VB\sigma}$}.
  \qquad
  \label{intsecondan}
\Esist
Now, we want to let $n$ tend to infinity in~\eqref{intsecondan}.
First, by \eqref{secondan} with $v=0$, we have~that
\Bsist
  && \norma{\Beta(\yn(t))}_{\Luno}
  \leq \norma{B^\sigma \yn(t)}^2
  + \iO \Beta(\yn(t))
  \non
  \\
  &&\leq \bigl(
    \norma{\tau\dt\yn(t)}
    + \norma{\pi(\yn(t))}
    + \norma{\un(t)}
    + \norma{\mun(t)}
  \bigr) \, \norma{\yn(t)}
  \quad \aat.
  \non
\Esist
So, by accounting for the \Lip\ continuity of $\pi$, and owing to \accorpa{stimeun}{stimemun},
we obtain \juerg{that}
\Beq
  \norma{\Beta(\yn)}_{\L2\Luno}
  \leq c_T \,.
  \label{stimaBeta}
\Eeq
On the other hand, by recalling \pier{\eqref{ynbdd}, \eqref{dtyntozero} and}
the compact embedding $\VB\sigma\subset H$ that follows \juerg{from} \eqref{hpAB}\pier{,
we can apply \cite[Sect.~8, Cor.~4]{Simon} and} deduce that 
\Beq
  \yn \to \yi
  \quad \hbox{strongly in $\C0H$}.
  \label{strongyn}
\Eeq
We infer that $\pi(\yn)$ converges to $\pi(\yi)$ in the same topology
since $\pi$ is \Lip\ continuous.
In order to deal with the nonlinearity~$\Beta$, we notice that
we can assume that $\yn\to\yi$ a.e.\ in $\Omega\times(0,T)$
so that, by lower semicontinuity, we deduce the inequality
\Beq
  \intQ \Beta(\yi)
  \leq \liminf_{n\to\infty} \intQ \Beta(\yn)
  \non
\Eeq
where the last term is finite by~\eqref{stimaBeta}.
As \eqref{convyn} also implies that
\Beq
  \ioT \norma{B^\sigma\yi(t)}^2 \, dt
  \leq \liminf_{n\to\infty} \ioT \norma{B^\sigma\yn(t)}^2 \, dt\,,
  \non
\Eeq
and \juerg{since the second statement in \eqref{stimeun} yields that} $\un\to\ui$ strongly in $\L2H$, \pier{from \eqref{intsecondan} and \eqref{convyn} it follows} that $\yi$ satisfies the variational inequality
\Bsist
  && \ioT \bigl( B^\sigma \yi(t) , B^\sigma( \yi(t)-v(t)) \bigr) \, dt 
  + \intQ \Beta(\yi) 
  \non
  \\
  && \quad {}
  + \ioT \bigl( \pi(\yi(t)) - \ui  ,  \yi(t)-v(t) \bigr)\,dt 
  \non
  \\
  && \leq \intQ \Beta(v)
  \qquad \hbox{for every $v\in\L2{\VB\sigma}$}.
  \label{intsecondai}
\Esist
Equivalently, $\yi$ \pier{fulfills}
\begin{align}
  & \bigl( B^\sigma \yi(t) , B^\sigma( \yi(t)-v) \bigr)
  + \iO \Beta(\yi(t))
  + \bigl( \pi(\yi(t)) \pier{{}-\ui{}} ,  \yi(t)-v \bigr)
  \leq \iO \Beta(v)
  \non
  \\
  & \quad \hbox{for every $v\in\VB\sigma$ and \aat}.
  \label{secondai}
\end{align}
At this point, we can easily conclude.
\pier{In view of} \eqref{dtyntozero} and \eqref{convyn},
we have that
\Beq
  \dt\yi = 0 \,,
  \quad \hbox{whence $\yi$ takes a constant value $\bar y\in\VB\sigma$ on $[0,T]$}.
  \non
\Eeq
On the other hand, $\yn(0)$ converges to $\yi(0)$ in~$H$ by~\eqref{strongyn}.
Thus, $\yn(0)$ converges to $\bar y$ in~$H$.
As $\yn(0)=y(\tn)$ converges weakly to $\yo$ in $\VB\sigma$ by assumption,
we conclude that $\bar y=\yo$, that is,
\Beq
  \yi(t) =  \yo
  \quad \hbox{for every $t\in[0,T]$}.
  \label{yiyo}
\Eeq
Therefore, \eqref{secondai} becomes~\eqref{yostat}.

\step Second part, second case

Assume now that $\lambda_1=0$.
Coming back to the proof just concluded,
we see that the assumption $\lambda_1>0$ has been used 
just to obtain~\eqref{stimemun}, its consequence~\eqref{stimaBeta},
and to make $\mun$ disappear in the limiting inequality~\eqref{yostat}.
Therefore, the same argument essentially applies in the case $\lambda_1=0$
(with the modifications that are needed to prove \eqref{charomega} instead of~\eqref{yostat}),
\juerg{provided}~we can derive a convergence property for~$\mun$
(in~place of~\eqref{stimemun}, which should be false now) and~\eqref{stimaBeta}.
To this end, we recall assumption \eqref{hpmz} and notice that 
it implies the existence of \juerg{some} \,$\delta>0$\, such that
\,$\mz\pm\delta$\, belong to~$D(\beta)$.
Then, as $\,v\,$ in~\eqref{secondan},
we choose the convex combination 
$\frac12\,(\mz\pm\delta)+\frac12\,\yn(t)$ 
(thus, with values in~$D(\Beta)$),
which gives $\yn(t)-v=\frac12(\yn(t)-\mz\mp\delta)$.
\pier{Thanks} to the convexity of~$\Beta$, we obtain \aat
\Bsist
  && \frac \tau 2 \, \bigl( \dt\yn(t) , \yn(t) - \mz \mp \delta \bigr)
  + \frac 12 \, \bigl( B^\sigma \yn(t) , B^\sigma( \yn(t) - \mz \mp \delta ) \bigr)
  \non
  \\
  && \quad {}
  + \iO \Beta(\yn(t))
  + \frac 12 \, \bigl( \pi(\yn(t)) - \un(t) , \yn(t) - \mz \mp \delta \bigr)
  \non
  \\
  && \leq \frac 12 \, \bigl( \mun(t) ,  \yn(t) - \mz \mp \delta \bigr)
  + \iO \Beta \bigl( \textstyle\frac12\,(\mz\pm\delta)+\frac12\,\yn(t) \bigr)
  \non
  \\
  && \leq \frac 12 \, \bigl( \mun(t) ,  \yn(t) - \mz \mp \delta \bigr)
  + \frac 12 \iO \Beta (\yn(t))
  + \frac 12 \iO \Beta(\mz\pm\delta).
  \non
\Esist
By multiplying by~$2$ and rearranging, we deduce that
\Bsist
  && \pm\delta \, \bigl( \mun(t),1 \bigr)
  + \norma{B^\sigma \yn(t)}^2
  + \iO \Beta(\yn(t))
  \non
  \\
  && \leq - \tau\, \bigl( \dt\yn(t) , \yn(t) - \mz \mp \delta \bigr)
  + \bigl( B^\sigma \yn(t) , B^\sigma( \mz \pm \delta ) \bigr)
  \non
  \\
  && \quad {}  
  - \bigl( \pi(\yn(t)) - \un(t) , \yn(t) - \mz \mp \delta \bigr)
  + \iO \Beta(\mz\pm\delta)
  \non
  \\
  && \quad {}  
  + \bigl( \mun(t) ,  \yn(t) - \mz \bigr).
  \label{perstimamu}
\Esist
Now, we recall the conservation property~\eqref{conservation} 
and \juerg{note that the Poincar\'e type inequality \eqref{poincare} is valid}  since $\lambda_1=0$.
We thus have that \aat \,\juerg{it holds}
\Bsist
  && \bigl( \mun(t) ,  \yn(t)-\mz \bigr)
  = \bigl( \mun(t) - \mean\mun(t) ,  \yn(t)-\mz \bigr)
  \nonumber
  \\
  && \leq c \, \norma{A^r(\mun(t)-\mean\mun(t))} \, \norma{\yn(t)-\mz}
  = c \, \norma{A^r\mun(t)} \, \norma{\yn(t)-\mz},
  \nonumber
\Esist
so that we can use \eqref{forsemstimeun} in the \rhs\ of~\eqref{perstimamu}.
By also accounting for \eqref{stimeun}, \eqref{ynbdd} and~\eqref{dtyntozerobis},
we deduce that the function
\Beq
  t \mapsto \delta \, \bigl| \bigl( \mun(t),1 \bigr) \bigr|
  + \iO \Beta(\yn(t))
  \non
\Eeq
is bounded in $L^2(0,T)$, uniformly with respect to~$n$.
In particular, \eqref{stimaBeta}~holds also in this case.
Moreover, the mean value of $\mun$ is estimated in $L^2(0,T)$
so that the definition \eqref{defnormaAr} of the norm in $\VA r$ and \eqref{forsemstimeun} imply that
$\mun$ is bounded in~$\L2{\VA r}$.
Therefore, we have\juerg{, at least for a subsequence,} 
\Beq
  \mun \to \mui
  \quad \hbox{weakly in $\L2{\VA r}$},
  \label{convmun}
\Eeq
which is the desired convergence property to be established in place of~\eqref{stimemun}.
At this point, we repeat the argument used in the case $\lambda_1>0$
provided that we modify~\eqref{secondai},
since we have \eqref{convmun} instead of~\eqref{stimemun}.
In place of that variational inequality, we obtain the following~one:
\Bsist
  && \bigl( B^\sigma \yi(t) , B^\sigma( \yi(t)-v) \bigr)
  + \iO \Beta(\yi(t))
  + \bigl( \pi(\yi(t)) - \ui ,  \yi(t)-v \bigr)
  \non
  \\
  && \leq \bigl( \mui(t) ,  \yi(t)-v \bigr)
  + \iO \Beta(v)
  \non
  \\
  && \quad \hbox{for every $v\in\VB\sigma$ and \aat}.
  \label{secondaibis}
\Esist
On the contrary, \eqref{yiyo}~holds true with the same proof also in the present case.
Finally, \eqref{forsemstimeun}~implies that $A^r\mui(t)=0$ \aat, i.e.,
that $\mui$ is space independent since $\lambda_1=0$ (\juerg{cf.}~\eqref{hpsimple}).
Therefore, \eqref{secondaibis} becomes
\Bsist
  && \bigl( B^\sigma \yo , B^\sigma (\yo-v) \bigr) 
  + \iO \Beta(\yo) 
  + \bigl( \pi(\yo) - \pier{\ui} , \yo-v \bigr)
  \non
  \\
  && \leq \bigl( \mui(t) , \yo-v \bigr)
  + \iO \Beta(v)
  \non
  \\
  && \quad \hbox{for every $v\in\VB\sigma$ and \aat}.
  \label{charomegaT}
\Esist
Hence, \pier{as \eqref{charomegaT} holds for arbitrary values of \,$T$, this}
implies \juerg{the validity of} \eqref{charomega} with a proper \juerg{function} 
$\mui\in L^\infty_{loc}([0,+\infty))$.
\pier{Indeed, let us} \juerg{denote for $\,m=0,1,\dots\,$} by $\,\mui^m\,$ the function $\,\mui\,$ 
\juerg{which satisfies} \eqref{charomegaT} with $\,T=m$,
and \juerg{we} construct $\,\mui\,$ on $\,(0,+\infty)\,$ by setting $\mui(t):=\mui^m(t)$
for a.a.\ $t\in(m,m+1)$, for $m=0,1,\dots$.
Then, \pier{it turns out that} $\mui\in L^\infty_{loc}([0,+\infty))$, and \pier{the inequality} \eqref{charomega} clearly holds.
This completes the proof of Theorem~\ref{Longtime}.\QED

\step
Proof of Proposition~\ref{Goodmui}

We come back to the proof just concluded\juerg{, keeping} its notation.
More precisely, we consider the part of the proof in the first case~$i)$ that also holds for the second one.
For the first part of Proposition~\ref{Goodmui}, 
it is sufficient to show that, for every $T\in(0,+\infty)$,
the function $\mui$ is unique and constant on $(0,T)$.
So, we fix an arbitrary $T>0$.
By virtue of the results of~\cite{CGS19} summarized just before the statement we are considering, 
we can replace the variational inequality \eqref{seconda} by the equation~\eqref{strong},
so that \eqref{secondan} can be written in the strong form
\Beq
  \tau \, \dt\yn + B^{2\sigma} \yn + \beta(\yn) + \pi(\yn) = \mun + \un
  \quad \hbox{a.e.\ in $\Omega\times(0,T)$}.
  \label{strongn}
\Eeq
Now, we observe that our assumption \eqref{gb} obviously implies that $\yn$ takes its values in~$[a,b]$.	
From \eqref{strongyn} and the \Lip\ continuity of $\beta+\pi$ in $[a,b]$, we thus infer that 
\Beq
  (\beta+\pi)(\yn) \to (\beta+\pi)(\yi)
  \quad \hbox{strongly in $\C0H$}.
  \non
\Eeq
On the other hand, by comparison in \eqref{strongn}, we see that $B^{2\sigma}\yn$ is in fact bounded in $\L2H$.
Therefore, the limiting function $\yi$ belongs to $\VB{2\sigma}$ and satisfies the equation
\Beq
  B^{2\sigma} \yi + \beta(\yi) + \pi(\yi) = \mui + \ui
  \quad \hbox{a.e.\ in $\Omega\times(0,T)$}.
  \non
\Eeq
But we already know that $\yi$ takes the constant value~$\yo$.
Therefore, $\yo\in\VB{2\sigma}$, and we have that
\Beq
  B^{2\sigma} \yo + \beta(\yo) + \pi(\yo) = \mui(t) + \ui
  \quad \hbox{for a.a.\ $t\in(0,T)$}.
  \non
\Eeq
By comparison, we conclude that $\mui$ is unique and time independent, thus constant, 
and the above equation becomes~\eqref{strongi}.
This completes the proof.

%%%%%%%%%%%%%%%%%%%%%%%%%%%%%%%%%%%%%%%%%%%%%%%%%%%%%%%%%%%%%%%%%%%%%%%%

\section*{Acknowledgments}
\pier{This research was supported by the Italian Ministry of Education, 
University and Research~(MIUR): Dipartimenti di Eccellenza Program (2018--2022) 
-- Dept.~of Mathematics ``F.~Casorati'', University of Pavia.}
PC and GG gratefully acknowledge some financial support from 
the GNAMPA (Gruppo Nazionale per l'Analisi Matematica, 
la Probabilit\`a e le loro Applicazioni) of INdAM (Isti\-tuto 
Nazionale di Alta Matematica) and the IMATI -- C.N.R. Pavia.

%%%%%%%%%%%%%%%%%%%%%%%%%%%%%%%%%
%% bibliography
%%%%%%%%%%%%%%%%%%%%%%%%%%%%%%%%%

\vspace{3truemm}

\Begin{thebibliography}{10}

\bibitem{AbWi}
H. Abels, M. Wilke,
Convergence to equilibrium for the Cahn--Hilliard equation with a logarithmic free energy,
{\em Nonlinear Anal.} {\bf 67} (2007), 3176-3193.

\bibitem{AkSS1}
G. Akagi, G. Schimperna, A. Segatti,
Fractional Cahn--Hilliard, Allen--Cahn and porous medium equations,
\pier{{\em J. Differential Equations}} {\bf 261} (2016), 2935-2985.

\bibitem{Barbu}
V. Barbu,
``Nonlinear Differential Equations of Monotone Type in Banach Spaces'',
Springer,
London, New York, 2010.

\bibitem{Brezis}
H. Brezis,
``Op\'erateurs maximaux monotones et semi-groupes de contractions
dans les espaces de Hilbert'',
North-Holland Math. Stud.
{\bf 5},
North-Holland,
Amsterdam,
1973.

\bibitem{ChFaPr}
R. Chill, E. Fa\v sangov\'a, J. Pr\"uss, 
Convergence to steady state of solutions of the Cahn--Hilliard and Caginalp equations with dynamic boundary conditions, 
{\em Math. Nachr.} {\bf 279} (2006), \pier{1448}-1462.

\bibitem{CG}
P. Colli, G. Gilardi,
Well-posedness, regularity and asymptotic analyses for a fractional phase field system,
\pier{{\em Asymptot. Anal.} \textbf{114} (2019), 93-128.}

\bibitem{CGLN}
P. Colli, G. Gilardi, \pier{Ph.} Lauren\c cot, A. Novick-Cohen,
Uniqueness and long-time behavior for the conserved phase-field system with memory,
{\em Discrete Contin. Dynam. Systems} {\bf 5} (1999), 375-390.

\bibitem{CGPS2}
P. Colli, G. Gilardi, P. Podio-Guidugli, J. Sprekels,
Well-posedness and long-time behaviour for a nonstandard viscous Cahn--Hilliard system,
{\em SIAM J. Appl. Math.} {\bf 71} (2011), 1849-1870.

\bibitem{CGS17}
P. Colli, G. Gilardi, J. Sprekels,
On the longtime behavior of a viscous Cahn--Hilliard system 
with convection and dynamic boundary conditions,
\pier{{\em J. Elliptic Parabol. Equ.} {\bf 4} (2018), 327-347.}

\bibitem{CGS18}
P. Colli, G. Gilardi, J. Sprekels,
Well-posedness and regularity for a generalized fractional Cahn--Hilliard system,
\pier{{\em Atti Accad. Naz. Lincei Rend. Lincei Mat. Appl.}, 
{\bf 30} (2019), 437-478.}

\bibitem{CGS19}
P. Colli, G. Gilardi, J. Sprekels,
Optimal distributed control of a generalized fractional Cahn--Hilliard system,
\pier{{\em Appl. Math. Optim.} doi:10.1007/s00245-018-9540-7 (see also 
preprint arXiv:1807.03218 [math.AP] (2018), pp.~1-36).}

\bibitem{EfGaZe}
M. Efendiev, H. Gajewski, S. Zelik, 
The finite dimensional attractor for a 4th order system of the Cahn--Hilliard type with a supercritical nonlinearity,
{\em Adv. Differential Equations} {\bf 7} (2002), 1073-1100.

\bibitem{EfMiZe}
M. Efendiev, A. Miranville, S. Zelik, 
Exponential attractors for a singularly perturbed Cahn--Hilliard system. 
{\em Math. Nachr.} {\bf 272} (2004), 11-31.

\bibitem{Gal1}
C.\,G. Gal,
Exponential attractors for a Cahn--Hilliard model in bounded domains with permeable walls,
{\em Electron. J. Differential Equations} \pier{(2006), No.~143, 23 pp.}

\bibitem{Gal3}
C. Gal,
Well-posedness and long time behavior of the non-isothermal viscous Cahn--Hilliard equation with dynamic boundary conditions,
{\em Dyn. Partial Differ. Equ.} {\bf 5} (2008), 39-67.

\bibitem{Gal2}
C.\,G. Gal,
Non-local Cahn--Hilliard equations with fractional dynamic boundary,
{\em European J. Appl. Math.} {\bf 28} (2017), 736-788.

\bibitem{GaGr1}
C.\,G. Gal, M. Grasselli,
Asymptotic behavior of a Cahn--Hilliard--Navier--Stokes system in~2D,
{\em Ann. Inst. H. Poincar\'e Anal. Non Lin\'eaire } {\bf 27} (2010), 401-436.

\bibitem{GaGr2}
C.\,G. Gal, M. Grasselli, 
Instability of two-phase flows: \pier{a} lower bound on the dimension of the global attractor  of the Cahn--Hilliard--Navier--Stokes system,
\pier{{\em Phys.~D}\/} {\bf 240} (2011), 629-635.

\bibitem{GiMiSchi}
\pier{G. Gilardi, A. Miranville, G. Schimperna,
Long-time behavior of the Cahn--Hilliard equation 
with irregular potentials and dynamic boundary conditions,
{\em Chin. Ann. Math. Ser. B\/} {\bf 31} (2010), 679-712.}

\bibitem{GiRo}
G. Gilardi, E. Rocca,
Well-posedness and long-time behaviour for a singular phase 
field system of conserved type,
{\em IMA J. Appl. Math.} {\bf 72} (2007), 498-530.

\bibitem{GS}
G. Gilardi, J. Sprekels,
Asymptotic limits and optimal control for the Cahn--Hilliard 
system with convection and dynamic boundary conditions,
{\em Nonlinear Anal.} {\bf 178} (2019), 1-31. 

\bibitem{GrPeSch}
M. Grasselli, H. Petzeltov\'a, G. Schimperna,  
Asymptotic  behavior  of  a  nonisothermal  viscous  Cahn--Hilliard  equation  with  inertial  term,
{\em J. Differential Equations} {\bf 239} (2007), 38-60.

\bibitem{JiWuZh}
J. Jiang, H. Wu, S. Zheng,
Well-posedness and long-time behavior of a non-autonomous Cahn--Hilliard--Darcy system with mass source modeling tumor growth,
{\em J. Differential Equations} {\bf 259} (2015), 3032-3077.

\bibitem{LiZho}
D. Li, C. Zhong,
Global attractor for the Cahn--Hilliard system with fast growing nonlinearity,
{\em J. Differential Equations} {\bf 149} (1998), 191-210.

\bibitem{Mir3}
A. Miranville, 
Long-time behavior of some models of Cahn--Hilliard equations in deformable continua, 
{\em Nonlinear Anal.} {\bf 2} (2001), 273-304.

\bibitem{Mir2}
A. Miranville, 
Asymptotic behavior of the Cahn--Hilliard--Oono equation, 
{\em J. Appl. Anal. Comput.} {\bf 1} (2011), 523-536.

\bibitem{Mir1}
A. Miranville, 
Asymptotic behavior of a generalized Cahn--Hilliard equation with a proliferation term, 
{\em Appl. Anal.} {\bf 92} (2013), 1308-1321.

\bibitem{MiZe2}
A. Miranville, S. Zelik, 
Robust exponential attractors for Cahn--Hilliard type equations with singular potentials, 
{\it Math. Methods Appl. Sci.} {\bf 27} (2004), 545-582.

\bibitem{MiZe1}
A. Miranville, S. Zelik, 
Exponential attractors for the Cahn--Hilliard equation with dynamic boundary conditions, 
{\em Math. Methods Appl. Sci.} {\bf 28} (2005), 709-735.

\bibitem{PrVeZa}
J.~Pr\"uss, V.~Vergara, R.~Zacher,
Well-posedness and long-time behaviour for the non-isothermal Cahn--Hilliard equation with memory,
\pier{{\em Discrete Contin. Dyn. Syst.}} {\bf 26} (2010), 625-647.

\bibitem{Seg}
A. Segatti, 
On the hyperbolic relaxation of the Cahn--Hilliard equation in~3D: approximation and long time behaviour,
{\em Math. Models Methods Appl. Sci.} {\bf 17} (2007), 411-437.

\bibitem{Simon}
J. Simon,
Compact sets in the space $L^p(0,T; B)$,
{\it Ann. Mat. Pura Appl.~(4)\/} 
{\bf 146} (1987), 65-96.

\bibitem{WaWu}
X.-M. Wang, H. Wu, 
Long-time behavior for the Hele--Shaw--Cahn--Hilliard system,
{\em Asymptot. Anal.} {\bf 78} (2012), 217-245.

\bibitem{WuZh}
H. Wu, S. Zheng, 
Convergence to equilibrium for the Cahn--Hilliard equation with dynamic boundary conditions, 
{\em J. Differential Equations} {\bf 204} (2004), 511-531.

\bibitem{Yeh}
C.-C. Yeh,
Discrete inequalities of the Gronwall-Bellman type in $n$ independent variables,
\emph{J. Math. Anal. Appl.} {\bf 105} (1985), 322-332.

\bibitem{ZhaoLiu}
X. Zhao, C. Liu,
On the existence of global attractor for 3D viscous Cahn--Hilliard equation,
{\em Acta Appl. Math.} {\bf 138} (2015), 199-212.

\bibitem{Zheng}
S. Zheng,
Asymptotic behavior of solution to the Cahn--Hilliard equation,
{\em Appl. Anal.} {\bf 23} (1986), 165-184.

\End{thebibliography}

\End{document}

%%%%%%%%%%%%%%%%%%%%%%%%%%%%%%%%%%%%%%%%%%%%%